\documentstyle[fullpage,12pt,psfig]{article}
\def\SBIMSMark#1#2#3{
 \font\SBF=cmss10 at 10 true pt
 \font\SBI=cmssi10 at 10 true pt
 \setbox0=\hbox{\SBF Stony Brook IMS Preprint \##1}
 \setbox2=\hbox to \wd0{\hfil \SBI #2}
 \setbox4=\hbox to \wd0{\hfil \SBI #3}
 \setbox6=\hbox to \wd0{\hss
             \vbox{\hsize=\wd0 \parskip=0pt \baselineskip=10 true pt
                   \copy0 \break%
                   \copy2 \break%
                   \copy4 \break}}
 \dimen0=\ht6   \advance\dimen0 by \vsize \advance\dimen0 by 8 true pt
                \advance\dimen0 by -\pagetotal
 \dimen2=\hsize \advance\dimen2 by .25 true in
  \openin2=publishd.tex
  \ifeof2\setbox0=\hbox to 0pt{}
  \else 
     \setbox0=\hbox to 3.1 true in{
                \vbox to \ht6{\hsize=3 true in \parskip=0pt  \noindent  
                \input publishd.tex 
                \vfill}}
  \fi
  \closein2
  \ht0=0pt \dp0=0pt
 \ht6=0pt \dp6=0pt
 \setbox8=\vbox to \dimen0{\vfill \hbox to \dimen2{\copy0 \hss \copy6}}
 \ht8=0pt \dp8=0pt \wd8=0pt
 \copy8
 \message{*** Stony Brook IMS Preprint #1, #2 ***}
}

\title{Ergodicity of conformal measures for unimodal polynomials}
\author{Eduardo A. Prado \thanks{Supported in part by CNPq-Brazil and 
S.U.N.Y. at Stony Brook}\\  
Instituto de Matem\'atica e Estat\'istica \\
Universidade de S\~ao Paulo \\Caixa Postal 66281 CEP 05389-970 \\
S\~ao Paulo, Brazil \\ e.mail: prado@ime.usp.br}

\date{ }

\newcommand{\finren}{ finitely many times renormalizable with only 
repelling periodic points}

\newcommand{\pollikef}{f: \bigcup U_{i} \rightarrow U}
\newcommand{\pollikeg}{g: \bigcup V_{i} \rightarrow V}

\newcommand{\seclimb}{\cal SL}

\newcommand{\yocpieceni}{Y^{n(i)}(x_{i})}
\newcommand{\yocpiecen}{Y^{n}(x)}

\newcommand{\yoccritpiecen}{Y^{n}(0)}
\newcommand{\yoccritpiecem}{Y^{j}(0)}
\newcommand{\yoccritpiecenmk}{Y^{n+k}(0)}

\newcommand{\yoccritpiecei}{Y^{i}(0)}
\newcommand{\yoccritpieceiu}{Y^{i_{1}}(0)}
\newcommand{\yocpiecenj}{Y^{n(j)}(x)}
\newcommand{\yocpieceyi}{Y^{i}(y)}

\newcommand{\julf}{ {\em J(f)}}
\newcommand{\quadrpol}{quadratic polynomial}

\newcommand{\renf}{R(f)}
\newcommand{\diam}{ {\rm diam} }
\newcommand{\mod} { {\rm mod} }

\newcommand{\yoccritdepn}{Y^{n}(0)} 
\newcommand{\yocdepni}{Y^{n}_{i}}

\newcommand{\yocdepzz}{V^{0,0}}
\newcommand{\yocdepzu}{V^{0,1}}
\newcommand{\yocdepzt}{V^{0,t}}
\newcommand{\yocdepztmu}{V^{0,t+1}}
\newcommand{\yocdepztz}{V^{0,t(0)}}
\newcommand{\yocdepztzmu}{V^{0,t(0)+1}}
\newcommand{\yocdepzn}{V^{0,n}}

\newcommand{\yocdepuz}{V^{1,0}}
\newcommand{\yocdepuu}{V^{1,1}}
\newcommand{\yocdeput}{V^{1,t}}
\newcommand{\yocdeputu}{V^{1,t(1)}}
\newcommand{\yocdeputumu}{V^{1,t(1)+1}}
\newcommand{\yocdepmz}{V^{m,0}}
\newcommand{\yocdepmu}{V^{m,1}}
\newcommand{\yocdepmtm}{V^{m,t(m)}}
\newcommand{\yocdepmtmmu}{V^{m,t(m)+1}}

\newcommand{\sulpiecejmu}{V^{n,1}}
\newcommand{\sulpiecej}{V^{n,0}}

\newtheorem{dummylemma}{Dummy}[section]

\newtheorem{proposition}[dummylemma]{Proposition}
\newtheorem{theorem}[dummylemma]{Theorem}
\newtheorem{thm}[enumi]{Theorem}
\newtheorem{lemma}[dummylemma]{Lemma}
\newtheorem{cor}[enumi]{Corollary}

\newtheorem{ExamplE}[dummylemma]{\em Example}
\newtheorem{DefinitioN}[dummylemma]{\em Definition}

\newenvironment{romanEnum}
                {\begin{list}
                    {(\roman{enumi})}           
                    {\usecounter{enumi}         
                     \setlength{\rightmargin}{\leftmargin}
                     \setlength{\labelwidth}{7.5mm}
                    }
                }{\end{list}}

\newcommand{\proof}[1]{{\noindent \bf Proof.~~}#1 \hspace{1cm} \medskip$\Box$\hfill}

\begin{document}

\maketitle
\SBIMSMark{1996/8}{June 1996}{}

\section{Introduction}
\label{intro}
Our goal in this work is to study ergodic properties of
polynomial dynamical systems using conformal
measures.

 Let $f: {\bf C} \rightarrow {\bf C} $ be a polynomial. Sullivan
showed in \cite{s} that it is possible to construct a conformal
measure for $f$ with support on $\julf$, the Julia set of $f$, for at
least one positive exponent $\delta$. By a conformal measure (or
$\delta$-conformal measure, to be more precise) we understand a Borel
probability measure $\mu$ satisfying the following condition:

\[\mu(f(A)) = \int_{A} |Df (z)|^{\delta} d\mu(z),\]
whenever $f$ restricted to the set $A$ is one to one.

We say that $\mu$ is ergodic if $\mu (X)=0$ or $\mu(X) =1$ whenever we
have $X=f^{-1}(X)$. Notice that usually when one talks about
ergodicity of a measure it is assumed that the measure is invariant.
In our case, due to the definition of conformal measure we are not
dealing with an invariant measure but rather a quasi-invariant
measure. We will show that a conformal measure is ergodic if $f$ is a
polynomial with certain properties. Let us introduce some definitions
in order to state our main result.

Let $f$ be a quadratic polynomial with only repelling periodic points.
Following \cite{l4}, we will say that $f$ satisfies the ${\it
secondary
\; limb\; condition}$ if there is a finite family of truncated
secondary limbs $L_{i}$ of the Mandelbrot set such that the hybrid
classes of all renormalizations $R^{m}(f)$ belong to $\cup L_{i}$. Let
$\seclimb$ stand for the class of quadratic polynomials satisfying
the secondary limb condition.

Some examples of polynomials of class $\seclimb$ are: Yoccoz and
Lyubich polynomials, and also infinitely many times renormalizable
real polynomials of degree two. A quadratic polynomial is a Yoccoz
polynomial if it is at most finitely many times renormalizable,
with only repelling periodic points. A Lyubich polynomial is an
infinitely many times renormalizable quadratic polynomial in
$\seclimb$ with sufficiently high combinatorics as described in
\cite{l}. 

If $f:U \rightarrow V$ is a polynomial-like map, then we say that
$\mod (f)$ is equal to the modulus of the topological annulus $A = V
\setminus U$. We say that an infinitely renormalizable polynomial 
$f$ has ${\it a \; priori \; bounds}$ if
there exists an $\epsilon >0$ such that $\mod (R^{m}(f)) > \epsilon$,
where $R^{m}(f)$ is the $m^{th}$ renormalization of $f$, for
infinitely many $m$. According to
\cite{gs}, \cite{lsa} and \cite{ly}, all infinitely
renormalizable real unimodal polynomials have a priori bounds.  By
{\em unimodal polynomial} we mean an even degree polynomial with just
one critical point.

We will show the following:

\begin{thm}
\label{ergthm}
Let $f$ be a polynomial with just one critical point 
and $J(f)$ connected.  Let $\mu$ be
a conformal measure for $f$. Suppose that $f$ is either:

\begin{romanEnum}
\item of class $\seclimb$, finitely many times renormalizable, or  
infinitely many times renormalizable with a priori bounds, or
\item a polynomial with a parabolic periodic
point, or
\item any unimodal polynomial with real coefficients.

\end{romanEnum}

\noindent Then $\mu$ is ergodic.
\end{thm}

\setcounter{enumi}{1}

Ergodicity of conformal measures is known if $f$ is expanding on
$\julf$ (see \cite{b},
\cite{s} and \cite{w}). Suppose that $f$ has just one critical point 
and an attracting cycle. Then $f$ is expanding on $J(f)$. If $J(f)$
is disconnected (and $f$ has just one critical point) then $f$ is again
expanding when restricted to $J(f)$. From this fact and
Theorem~\ref{ergthm} we conclude the following:

\begin{cor}
If $f$ is any unimodal polynomial with real coefficients
and $\mu$ a conformal measure for $f$, then $\mu$ is ergodic.
\end{cor}

The following is an immediate consequence of ergodicity of
conformal measures:

\begin{cor}
\label{uniq}
Let $f$ be as in Theorem~\ref{ergthm}. Then for any $\delta >0$,
there exists at most one $\delta$-conformal measure for $f$.
\end{cor}

The situation studied in this paper is the complex counterpart for the
ergodicity result in \cite{bl} where the Lebesgue measure is shown to
be ergodic under $S$-unimodal maps.\\

\noindent {\bf Acknowledgment.} I would like to thank Misha Lyubich for many
hours of valuable mathematical conversations and for continuous
encouragement. I also thank Marco Martens for various helpful
discussions. I thank John Milnor, Feliks Przytycki and Mariuzs
Urba\'{n}ski for reading a preliminary version of this paper and
making important comments.  I am grateful to CNPq-Brazil and the
Department of Mathematics of SUNY at Stony Brook for financial
support. Part of this work its taken from the authors's PhD thesis
presented in June 1995 at SUNY at Stony Brook under the direction of
M.  Lyubich.

\section{Renormalization and combinatorics}
\label{yocpyz}

\subsection{Non-renormalizable polynomials}
\label{yocpol}

We will briefly describe how to construct the Yoccoz puzzle pieces for
a quadratic polynomial. See \cite{h} and \cite{m} for a complete
exposition of such construction.

In this section we will consider quadratic polynomials $f$ with
repelling periodic points. We shall keep in mind though that the
construction of Yoccoz puzzles that will be described in this section 
can be repeated for polynomials with degree greater than two.

We say that $g:U \rightarrow U'$ is a ${\em quadratic-like \: map}$ if
it is a double branched covering and $U$ and $U'$ are open topological
discs with $U$ compactly contained in $U'$. The ${\em filled \:in \:
Julia \: set}$ of $g$ is the set $\{ z \in U$ : $g^{n}(z) $ is defined
for all natural numbers $n\}$. There are two fixed points of $g$
inside its filled in Julia set. If the filled in Julia set is
connected and both fixed points are repelling, one of them, the
dividing fixed point, disconnects the filled in Julia set of $g$ in
more than one connected component. The other does not.  Usually the
dividing fixed point is denoted by $\alpha$.  Quadratic-like maps were
first introduced and studied in \cite{dh}.

Remember that a polynomial or a polynomial-like map $f$ with connected
filled in Julia set is ${\em renormalizable}$ if there exist open
topological discs $U \subset U'$ with $0 \in U$ with $R(f):U
\rightarrow U'$ being a quadratic-like map with connected filled in
Julia set. We define $R(f)=f^{k}|_{U} $, with $k$ the smallest natural
number bigger than 1 satisfying the previous conditions (we call $k$
the {\em period of renormalization}).  Here $R(f)$ stands for the
renormalization of $f$. We can ask whether $R(f)$ is renormalizable or
not and then define renormalizations of $f$ of higher orders. So, each
renormalization of $f$ defines a quadratic polynomial-like map. We
refer the reader to \cite{mc} for more details concerning
renormalization.

Let $f$ be a degree two non-renormalizable polynomial with both fixed
points repelling and let $G$ be the Green function of the filled Julia
set of $f$. There are $q$ external rays landing at the dividing fixed
point of $f$, where $q
\geq 2$. The $q$ {\em Yoccoz puzzle pieces of depth zero} are the 
components of the topological disc defined by $G(z) <G_{0}$, where
$G_{0}$ is any fixed positive constant, cut along the $q$ external
rays landing at the dividing fixed points. We denote $Y^{0}(x)$ the
puzzle piece of depth zero containing $x$. We define the {\em puzzle
pieces of depth n} as being the connected components of the pre-images
of any puzzle piece of depth zero under $f^{n}$. Again, if $x$ is an
element of a given puzzle piece of depth $n$ we denote such puzzle
piece by $Y^{n}(x)$.

Suppose now that $f$ is at most finitely renormalizable with only
repelling periodic points. Let $\alpha$ be the dividing fixed point
of the last renormalization of $f$. Let $G$ be the Green function of
the filled in Julia set of $f$. In that case we define the puzzle
pieces of depth zero as being the components of the topological disc
$G(z) < G_{0}$, $G_{0}$ a positive constant, cut along the rays
landing at all points of the $f$-periodic orbit of $\alpha$. As before
we define the puzzle pieces of depth $n$ as being the connected
components of the pre-images under $f^{n}$ of the puzzle pieces of
depth zero. The puzzle piece at depth $n$ containing $x$ is denoted by
$Y^{n}(x)$.

We will consider the Yoccoz puzzle pieces as open topological discs.
Under this consideration the Yoccoz partition will be well defined
over the Julia set of the polynomial $f$ minus the set of pre-images
of the dividing fixed point of the last renormalization of $f$ (which
is $f$ itself in the non-renormalizable case).

A quadratic polynomial is a ${\em Yoccoz\: polynomial}$ if it is at
most finitely renormalizable with only repelling periodic points. We
will need the following result:

\begin{theorem}[Yoccoz]
\label{yocthm}
If f is a Yoccoz polynomial then $\bigcap_{n \geq 0} Y^{n}(x) = \{ x\}$
for any x where the Yoccoz partition is defined.
\end{theorem}

The following is the analog Theorem for higher degree real unimodal
polynomials:

\begin{theorem}[\cite{lsb}]
\label{lsthm}
Let $f(z)=z^{l} +c$, $l$ even, $c$ real and $f$ \finren. Then for any
$x$ where the Yoccoz partition is defined we have: $\bigcap_{n \geq 0}
Y^{n}(x) = \{ x\}$.
\end{theorem}

\subsection{The $\seclimb$ class}
\label{lyupol}

Here we will describe the secondary limb class of quadratic polynomials.
See \cite{l3} and \cite{l4} for a detailed exposition on this matter.
We will need some technical definitions.

Let us start with a quadratic polynomial $f$ with only repelling
periodic points. Given a Yoccoz puzzle piece $\yocdepni$ of $f$ and a
point $x$ such that $f^{j}(x)$ belongs to $\yocdepni$. We define the
{\em pull back of} $\yocdepni$ {\em along the orbit of x} as being the
only connected component of $f^{-j}(\yocdepni )$ containing $x$.  If
moreover $x$ belongs to $\yocdepni $ and $j$ is minimal and non-zero,
then we say that $j$ is {\em the first return time of x to}
$\yocdepni$. A puzzle piece is said to be a {\em critical puzzle
piece} if it contains the critical point. Notice that if we pull back
a critical puzzle piece $\yoccritdepn$ along the first return of the
critical point to $\yoccritdepn$ we get a new critical puzzle piece.

Suppose that $f$ is not Douady-Hubbard immediately renormalizable (see
\cite{l3}). Then it is possible to find a critical puzzle 
piece (that will be denoted by $\yocdepzz$) satisfying the following:
if the pull back of $\yocdepzz$ along the first return of the critical
point to $\yocdepzz$ is denoted by $\yocdepzu$, then the closure of
$\yocdepzu$ is properly contained in $\yocdepzz$. We keep repeating
this procedure: define $\yocdepztmu$, the puzzle piece of level $t+1$,
as being the pull back of $\yocdepzt$, the puzzle piece of level $t$,
along the first return of the critical point to $\yocdepzt$.  This
procedure stops if the critical point does not return to a certain
critical puzzle piece. If we assume that the critical point is
combinatorially recurrent, then we can repeat this procedure forever.
So let us assume that case. The collection $\yocdepzt$ for $t$ being a
natural number is the {\em principal nest of the first renormalization
level}.

Now we have a sequence of first return maps $f^{l(t)}:\yocdepztmu
\rightarrow \yocdepzt$. By definition $\yocdepzz$ properly contains
$\yocdepzu$. This implies that each $\yocdepzt$ properly contains
$\yocdepztmu$. It is also easy to see that each $f^{l(t)}:\yocdepztmu
\rightarrow \yocdepzt$ is a quadratic-like map.

We say that $f^{l(t)}:\yocdepztmu \rightarrow \yocdepzt$ is a ${\em
central\: return}$ or that $t$ is a ${\em central\: return\: level}$
if $f^{l(t)}(0)$ belongs to $\yocdepztmu$. A ${\em cascade\: of\:
central\: returns}$ is a set of subsequent central return levels. More
precisely, a cascade of central returns is a collection of central
return levels $t=t_{0}, ...,t_{0}+(N-1)$ followed by a non-central
return at level $t_{0}+N$. In this case we say that the above
cascade of central returns has length $N$. We could also have an
infinite cascade of central returns.  Notice that with the above
terminology a non-central return level is a cascade of central return
of length zero.

It is possible to show that the principal nest of the first
renormalization level ends with an infinite cascade of central returns
if and only if $f$ is renormalizable (see \cite{l3}). In that case,
denote the first level of this infinite cascade of central returns by
$t(0)$. Then we define the first renormalization $\renf$ of $f$ as
being the quadratic-like map $f^{l(t(0))}:\yocdepztzmu
\rightarrow \yocdepztz$. The filled-in Julia set of $\renf$ is 
connected (it is also possible to show that $\bigcap \yocdepzn=
J(\renf)$). Again we can find the dividing fixed point of the Julia
set of $\renf$, some external rays landing at it and define new puzzle
pieces over the Julia set of $\renf$. The rays landing at the new
dividing fixed point are not canonically defined (remember that
$\renf$ is a polynomial-like map). We are not taking the external rays
of the original polynomial. Instead we need to make a proper selection
of those rays (see
\cite{l3}). As before we can construct the principal nest for $\renf$,
provided that $\renf$ is not Douady-Hubbard immediately
renormalizable. The elements of this new principal nest are denoted
by $\yocdepuz, \yocdepuu, ..., \yocdeput, ...$ and the nest is called
the {\em principal nest of the second renormalization level}. If this
new principal nest also ends in an infinite cascade of central
returns, we repeat the procedure just described and construct a third
principal nest. We repeat this process as many times as we can.

Now we define the {\em principal nest} of the polynomial $f$ as
being the set of critical puzzle pieces

\[\yocdepzz \supset \yocdepzu \supset ...\supset \yocdepztz \supset 
\yocdepztzmu \supset \yocdepuz \supset \yocdepuu \supset  ..., \yocdeputu 
\supset \]

\[ \supset \yocdeputumu \supset ... \supset \yocdepmz \supset \yocdepmu 
\supset ...\supset 
\yocdepmtm \supset \yocdepmtmmu \supset ...\]

In order to go ahead with the definition of the class of polynomials
we are interested in, we need the notion of a truncated secondary
limb.  A {\em limb} in the Mandelbrot set $M$ is the connected
component of $M \setminus \{ c_{0} \}$ not containing $0$, where
$c_{0}$ is a bifurcation point on the main cardioid. If we remove from
the limb a neighborhood of its root $c_{0}$, we get a {\em truncated
limb}. A similar object corresponding to the second bifurcation from
the main cardioid is a {\em truncated secondary limb}.

We say that a quadratic polynomial with only repelling periodic
points satisfies {\em the secondary limb condition} if there is a
finite family of truncated secondary limbs $L_{i}$ of the Mandelbrot
set such that the hybrid class of all renormalizations $R^{m}(f)$
belongs to $\cup L_{i}$. Let $\seclimb$ stand for the class of
quadratic polynomials satisfying the secondary limb condition.

\subsection{Complex Bounds}

We say that $f$ has {\em a priori bounds or complex bounds} if there
exists an $\epsilon >0$ such that $\mod (R^{m}(f)) >
\epsilon$, for infinitely many renormalizations of $f$.  A priori
bounds is, as we shall see, one of the main properties that we will be
using in this paper. In \cite{l4} it was conjectured that the
secondary limb condition described above implies a priori bounds. In
\cite{l} and \cite{l4} it was constructed a large class of infinitely many
times renormalizable quadratic polynomials satisfying the secondary
limb condition with a priori bounds. The next Theorem follows from
Theorem II in \cite{l3}.

\begin{theorem}[Lyubich]
\label{Lyuthm}
Let f be a Yoccoz polynomial. There exists a constant $c >0$
such that $\mod (\yocdepzt \setminus \yocdepztmu  )> c$, for all $t$.
\end{theorem}

Complex bounds were proved first by Sullivan (see \cite{ms}) for
infinitely renormalizable quadratic maps with bounded combinatorics. In
\cite{gs}, \cite{lsa} and \cite{ly} the restriction on the combinatorics 
was removed. Also \cite{lsa} provides complex bounds for infinitely
renormalizable polynomials of the form $f(z)=z^{l}+c$, where $l$ is
even and $c$ is real.

\begin{theorem}
\label{complexboundsthm}
Let $f(z) = z^{l} + c$ be an infinitely renormalizable real polynomial
of even degree $l$. If $a_{n}$ is the period of the $n^{th}$
renormalization of $f$, then there exist topological discs
$\sulpiecej$ and $\sulpiecejmu$ such that:

\begin{romanEnum}
\item $0 \in \sulpiecejmu$;
\item $ {\rm cl}(\sulpiecejmu )\subset \sulpiecej$;
\item $\mod( \sulpiecej \setminus \sulpiecejmu) \geq c > 0$;
\item $ f^{a_{n}}: \sulpiecejmu \rightarrow \sulpiecej$ is a polynomial-like 
map of degree $l$ with connected filled in Julia set;
\item $\diam(\sulpiecej) \rightarrow 0$ as $ n \rightarrow \infty$.
\end{romanEnum}
\end{theorem}

\subsection{Unbranched maps}
\label{unbranchedrenorm}

\begin{DefinitioN}[\cite{l1}]
\label{pollikedef}
Let $U$ and $U_{i}$ be open topological discs, $i=0, 1, ...,n$.
Suppose that $cl(U_{i}) \subset U$ and $U_{i} \bigcap U_{j} =
\emptyset$ if $i$ is different than $j$. A generalized polynomial-like
map is a map $\pollikef$ such that the restriction $f|U_{i}$ is a
branched covering of degree $d_{i}, d_{i} \geq 1$.
\end{DefinitioN}

We will not use the above Definition in full generality. From now on,
all generalized polynomial-like maps in this work will have just one
critical point. We will fix our notation as follows: $f|U_{0}$ is a
branched covering of degree $d$ onto $U$ (with zero being the only critical
point) and $f|U_{i}$ is an isomorphism onto $U$, if $i=1,...,n$.

For the next definition we will consider a polynomial $f$ with just one
critical point and $\pollikeg$ a generalized polynomial-like map.  Assume that
$g$ is defined in each $V_{i}$ as the first return map to $V$ under $f$ of the
elements of $V_{i}$. In particular we assume that $V_{0}$ is the pull back of
$V$ along the first return of the critical point of $f$ to $V$.

\begin{DefinitioN}[\cite{lsb} and \cite{mc}]
We say that $g$ is unbranched if whenever $f^{i}(0)$ belongs to $V$,
then $f^{i}(0)$ is an iterate of $0$ under $g$.
\end{DefinitioN}

Let $g$ be as in the above definition. If we assume that the critical
point is recurrent then the intersection of the critical set of $f$
with $V$ is contained in the domain of $g$. Notice that if $f$ is
renormalizable and if $g$ is a renormalization of $f$, then the above
definition coincides with the unbranched renormalization definition
from
\cite{mc}.

\begin{lemma}
\label{unbranchedprop}
Suppose that $\pollikeg$ is unbranched. Let $z$ be an element in $J(f)$
and $k$ be the smallest positive number such that $f^{k}(z)$ belongs to
$V_{0}$. Then we can pull $V$ back univalently along the orbit $z, f(z),
..., f^{k}(z)$.
\end{lemma}
\proof{ Suppose this is not the case. Then for some positive $r$ smaller than
$k$, $f^{-r}(V)$ hits the critical value $f(0)$. By $f^{-r}(V)$ we understand
the pull back of $V$ along the orbit $f^{k-r}(z), ..., f^{k}(z)$.  In that
case, $f^{i+1}(0) \in f^{i}(f^{-r}(V))$, for $i=0,..., r$. Suppose that there
exists $0 \leq i_{max} < r$ maximal such that $f^{i_{max} +1}(0)$ belongs to
$V$. In that case, by the unbranched property there exists a component $V_{j}$
of the domain of $g$ containing $f^{i_{max} +1}(0)$. As $V_{j}$ is the pull
back of $V$ under the first return of $f^{i_{max} +1}(0)$ to $V$, we conclude
that $f^{i_{max}}(f^{-r}(V)) \subset V_{j}$. Now, as the pull back of $V_{j}$
under $f$ along $0, ..., f^{i_{max}}(0)$ is contained in $V_{0}$, it follows
that $f^{-1}(f^{-r}(V)) \subset V_{0}$, contradicting the minimality of $k$.
Suppose now that there is no $i_{max}$, i.e., if the first return time of $0$
to $V$ is $r+1$. Then $f^{r}:f^{-r}(V) \rightarrow V$ is the univalent branch
of the first return map $f^{r+1} : V_{0} \rightarrow V$. This also contradicts
the minimality of $k$.  }

\begin{lemma}
\label{unbranchedlemma}
Let $f$ be an $\seclimb$ polynomial with a priori bounds. Then for
infinitely many $n$ we can find $U_{n}$ and $V_{n}$ such that the
$n^{th}$-renormalization of $f$ is given by $R^{n}(f): U_{n}
\rightarrow V_{n}$ and both the unbranched condition and the a priori
bound condition is verified for $U_{n}$ and $V_{n}$
\end{lemma}

\proof{
See Lemma 2.3 in \cite{l4} for the proof. 
}

Let $ f^{a_{n}}: \sulpiecejmu \rightarrow \sulpiecej$ be the polynomial-like 
maps introduced on Theorem~\ref{complexboundsthm} (the $n^{th}$
renormalization of $f(z)=z^{l}+c$). We have the following:

\begin{lemma}
\label{unbranchedinfrenlemma}
Let $f(z)= z^{l}+c$, with $l$ even and $c$ real. Then the polynomial-like maps
$ f^{a_{n}}: \sulpiecejmu \rightarrow \sulpiecej$ are unbranched, for
infinitely many $n$.
\end{lemma}

\proof{
This is due to the construction of the set $\sulpiecejmu$ and
$\sulpiecej$ in \cite{lsa}.
}

For the next Lemma, let $f$ be either a Yoccoz polynomial or any
finitely many times renormalizable real polynomial with only repelling
periodic points and just one critical point. Suppose that this
critical point is recurrent.  Let $\yoccritpiecen$ be a critical
Yoccoz piece and $\yoccritpiecenmk$ be the pull back of
$\yoccritpiecen$ along of the first return of the critical point to
$\yoccritpiecen$.

\begin{lemma}
\label{yocpullbacklemma}
Let be $z$ an element of $J(f)$ and let $m$ be the
smallest time that $ f^{m} (z) = \yoccritpiecenmk$. Then we can
univalently pull $\yoccritpiecen$ back along the orbit $z ,\ldots ,f^m
(z)$.
\end{lemma}

\proof{
If not, $f^{-t}(\yoccritpiecen )$ would contain the critical point, for
some $t$ less than $m$ (here $f^{-t}$ means the branch of $f^{-t}$
along the orbit of $x$). That would mean that $t$ is greater or equal
to the first return time of 0 to $\yoccritpiecen $. That would imply
$f^{-t}(\yoccritpiecen) \subset \yoccritpiecenmk$ by the Markov property of
puzzle pieces.  In other words, $z$ would hit $\yoccritpiecenmk$ on a time
strictly less than $m$, contradicting the definition of $m$.  
}

\begin{lemma}
\label{unbranchednonrenlemma}
Let $f$ be either a Yoccoz polynomial or any finitely many times
renormalizable real polynomial with only repelling periodic points.
Then for any critical puzzles $\yoccritpiecem$, we can find a
topological disk $D_{j}$ such that:

\begin{enumerate}
\item $\yoccritpiecem \subset D_{j}$
\item $\mod (D_{j} \setminus \yoccritpiecem ) > c(f)>0$
\item if $z \in J(f)$ and $m$ be the
smallest time that $ f^{m} (z) = \yoccritpiecem$. Then we can
univalently pull $D_{j}$ back along the orbit $z ,\ldots ,f^m
(z)$. 
\end{enumerate}
\end{lemma}

\proof{
For a Yoccoz polynomial this follows from Lemma~\ref{yocpullbacklemma}
and Theorem~\ref{Lyuthm}. If $f$ is real but not of degree two, then
the Lemma follows from the proof of Theorem B in \cite{lsa} (if the critical
orbit is minimal), and from Proposition 1.2 in \cite{lsb} if the
critical orbit is not minimal.  }

\section{Density Estimates}
\label{densest}

From now on $f$ will be a polynomial of even degree with just one
critical point and $\mu$ will denote a $\delta-$conformal measure
concentrated on the Julia set of $f$.

The analytic tool that we will use is the well known Koebe
distortion Theorem:

\begin{theorem}[Koebe]
\label{koebethm}
Let $A \subset B$ be two topological 
discs contained in the complex plane. Suppose that f is univalent when
restricted to B. Also suppose that \mbox{$B \setminus A$} is a
topological annulus with positive modulus m. Then

\[\frac{1}{K} \leq \frac{|Df(z_{1})|}{|Df(z_{2})|} \leq K\]
for all $z_{1}$ and $z_{2}$ in A, where the constant K depends
only on the number m.
\end{theorem}

The constant $K$ that appears in the Lemma is called the Koebe
constant.  Under the conditions of the above Lemma we say that $f$ has
bounded distortion inside the set $A$.

Let $f$ be either a Yoccoz polynomial or a finitely many times
renormalizable real polynomial with only repelling periodic points.
Notice that if a periodic point of $f$ in $\julf$ is expanding then
the set of all its pre-images has zero $\mu$-measure.  As we used just
expanding periodic points to construct puzzle pieces, given any closed
subset $X$ of $\julf$, we can create a cover $K_{i}$ of $X$ (up to a
set of zero measure) built up by puzzles pieces and with ${\rm lim}
\mu (K_{i}) =
\mu(X)$. This follows from Theorem~\ref{yocthm}, Theorem~\ref{lsthm} and 
the regularity of conformal measures.

\begin{DefinitioN}
The density of a set $X$ inside a set $Y$ is defined as follows:
${\rm dens}( X | Y ) = \frac{ \mu (X \bigcap Y)}{ \mu (Y)}$.
\end{DefinitioN}

\begin{lemma}
\label{existdensptlemma}
Let $f$ be either a Yoccoz polynomial or a finitely many times
renormalizable real polynomial with only repelling periodic points.
Let $X \subset J(f)$ be any measurable subset. If $\mu(X) >0$, there
is $x$ in $X$ such that ${\rm limsup} ({\rm dens} (X | \yocpiecen))
=1$.\end{lemma}

\proof{
Assume $\mu (X) >0$. If $X$ is not closed, take $W \subset X$ compact
with $\mu(X \setminus W)$ small. Notice that ${\rm dens} (X |
\yocpiecen) \geq {\rm dens} (W | \yocpiecen )$ for any $\yocpiecen$.
For all $\varepsilon >0$, there exists $i(\varepsilon )$, such that $1-
\varepsilon \leq \frac{ \mu (W \bigcap K_{i})}{ \mu (K_{i})} \leq 1$
if $i> i(\varepsilon)$ (remember that $K_{i}$ are the covers of $X$
made out of puzzle pieces). So we have for $i$ big ${\rm dens}(W |
K_{i} )= \frac{ \mu (W \bigcap K_{i})}{
\mu (K_{i})} \geq 1-
\varepsilon$.  As $K_{i}$ is the union of puzzle pieces we can
certainly find a puzzle piece in $K_{i}$, say $\yocpieceni$ such that
${\rm dens}(W | \yocpieceni) \geq 1-\varepsilon$.  Now replacing $X$
by $X \bigcap \yocpieceni$ and repeating this argument we will end up
with the desired result.  
}

\begin{DefinitioN}
The point $x \in X$ obtained in the previous Lemma is called a weak
density point of $X$.
\end{DefinitioN}

\begin{proposition}
\label{denstransfprop}
Let $A \subset B$ be two $\mu$-measurable subsets of the complex
plane. Suppose that f restricted to an open neighborhood of B is one
to one. Also suppose that there exists a positive constant K such
that

\[ \frac{1}{K} \leq \frac{|Df(z_{1})|}{|Df(z_{2})|} \leq K\]
for all $z_{1}$ and $z_{2}$ in B, then

\[ \frac{1}{K^ \delta} {\rm dens}(A|B) \leq {\rm dens}(f(A)|f(B)) \leq
K^ \delta {\rm dens} (A|B). \]

\end{proposition}

\proof{
Follows from the definition of conformal measure and the definition of
${\rm dens}(A | B)$.  }

If $U$ is a subset of the complex plane, we will denote by $U^{c}$ the
complement of $U$ inside the complex plane.

\begin{lemma}
\label{zeromeasyoclemma}
Let f be either a Yoccoz polynomial or a finitely many times
renormailzable real polynomial of even degree with only repelling
periodic points. Let $\mu$ a conformal measure for f.  Let U be any
neighborhood of the critical point. Then the set

\[ \{ x \in {\bf C} : f^{n} (x) \in U^{c}, \: for \: 
all \: positive \: n \} \]
has zero $\mu$-measure.
\end{lemma}

\proof{
It is enough to show this Lemma for $U= \yoccritpiecei$ because by
Theorem~\ref{yocthm} and Theorem~\ref{lsthm} any neighborhood of the
critical point contains some $\yoccritpiecei$, for $i$ sufficiently
big.  Suppose that the set $A =
\{ x \in {\bf C}$ : $f^{n} (x) \in \yoccritpiecei ^{c}$, for all $n$
positive $\}$ has positive measure, for some $i$ fixed. Then this set
has a point of weak density $x$, according to
Lemma~\ref{existdensptlemma}.  So we can find some sequence $n(j)
\rightarrow \infty$ such that ${\rm dens}(A |
\yocpiecenj ) \rightarrow 1$.

Notice that $f^{n(j)-i} ( \yocpiecenj )$ is a puzzle piece of depth
$i$ and none of the puzzle pieces $\yocpiecenj, f( \yocpiecenj )
\ldots , f^{n(j)-i} ( \yocpiecenj )$ contains the critical point. That
is because of the Markov property of puzzle pieces and the fact that
$\yocpiecenj$ contains elements of the set $A$. So for all
$\yocpiecenj$, $f^{n(j)-i}( \yocpiecenj)$ is a puzzle piece of depth
$i$ distinct from $\yoccritpiecei$ and the restriction $f^{n(j)-i}:
\yocpiecenj
\rightarrow Y^{i}(f^{n(j)-i}(x))$ is an isomorphism. As there exist 
just finitely many puzzle pieces of depth $i$ then there is a fixed
puzzle piece $\yocpieceyi$ (distinct from the one containing the
critical point) such that $f^{n(j)-i}(\yocpiecenj)= \yocpieceyi$ for
infinitely many $n(j)$. Passing to a subsequence and keeping the same
notation we will assume that $f^{n(j)-i}: \yocpiecenj\rightarrow
Y^{i}(y)$ is an isomorphism for all $n(j)$.

We will construct a neighborhood of $\yocpieceyi$ where the inverse
branch $f^{-(n(j)-i)}$ along the orbit $x, f(x), \ldots , f^{n_(j)
-i}(x)$ is defined as an isomorphism.

Let $i_{1}>i$ such that ${\rm mod}(\yoccritpiecei \setminus
\yoccritpieceiu)$ is positive. This is possible by Yoccoz's 
Theorem.

The boundary of $\yocpieceyi$ is composed by pairs of external rays
landing at points in the Julia set and equipotentials. The
intersection of this boundary with the Julia set is finite. Let $z$ be
a point of such finite intersection. Consider all puzzle pieces of
depth $i_{1}$ containing $z$ on its boundary. The closure of the union
of those puzzle pieces is a neighborhood of $z$ in the plane. Let us
call such neighborhood $V_{z}$. Notice that each equipotential and the
pieces of external rays landing at $z$ outside $V_{z}$ are at some
definite distance from the Julia set.  Take a small tubular
neighborhood (not intersecting the Julia set) of each one of the
equipotentials and pieces of external rays contained in the boundary
of $\yocpieceyi$, but outside $V_{z}$.  Now we define the neighborhood
$N$ of $\yocpieceyi$ as being the union of each $V_{z}$ with all
tubular neighborhoods described above and $\yocpieceyi$ itself (see
Figure~\ref{neighborhood}). Notice that we can make $N$ into a
topological disc if $i_{1}$ is big and the tubular neighborhoods
small. Also notice that since the distance between the boundaries of
$\yocpieceyi$ and $N$ is strictly positive, we get that ${\rm mod}(N
\setminus \yocpieceyi)$ is strictly positive.

\begin{figure}[bth]
\centerline{\psfig{figure=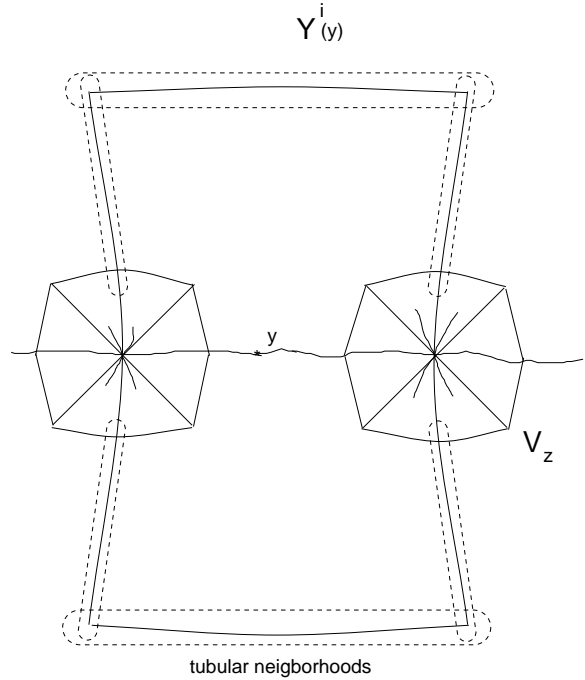,height=9cm}}
\caption{Construction of the neighborhood $N$ of $\yocpieceyi$}
\label{neighborhood}
\end{figure}

Now let us prove that we can pull $N$ back isomorphically along the
orbit $x, \ldots , f^{n(j)-i}(x)$ for any $n(j)$.

The pull back of $\yocpieceyi$ along $x, \ldots , f^{n(j)-i}(x)$
cannot hit the critical point. This is because $\yocpiecenj$ contains
points in the set $A$, the $f$-invariance of the set $A$ and the
Markov property of puzzle pieces. None of the pull backs of the
tubular neighborhoods can hit the critical point because those
neighborhoods are outside the Julia set. The pull back of $V_{z}$
along $x, \ldots , f^{n(j)-i}(x)$ can not touch the critical point.
If the pull back of $V_{z}$ would hit the critical point, then it
would intersect the interior of $\yoccritpieceiu$ (because $V_{z}$ is
made out of puzzle pieces of depth $i_{1}$). By the choice of $i_{1}$
(${\rm mod}(
\yoccritpiecei 
\setminus \yoccritpieceiu)>0$) 
 and because all the puzzle pieces of $V_{z}$ have a common boundary point
 with $\yocpieceyi$ we would conclude that some pre-image of
$\yocpieceyi$ along $x, \ldots , f^{n(j)-i}(x)$ would intersect $\yoccritpiecei$.
Contradiction!

So we can pull $N$ back isomorphically along the orbit $x, \ldots ,
f^{n(j)-i}(x)$ for any $n(j)$. By the construction of $N$ we have:
${\rm mod}(N \setminus \yocpieceyi)>0$. So we conclude that
$f^{n(j)-i}: \yocpiecenj \rightarrow
\yocpieceyi$ has bounded distortion with the Koebe constant not
depending on $n(j)$.

Using the above bounded distortion property,
Proposition~\ref{denstransfprop} and the fact that $x$ is a density
point for $A$, we conclude that ${\rm dens}(A| \yocpieceyi)$ is
arbitrarily close to one. On the other hand there exists some
pre-image of $\yoccritpiecei$ inside $\yocpieceyi$, so ${\rm dens}(A|
\yocpieceyi)$ is bounded away from 1. Contradiction!  }

Let us prove a similar result for the classes of infinitely
renormalizable polynomials that we are dealing with:

\begin{lemma}
\label{zeromeasinfrenormlemma}
Let f be any $\seclimb$
\quadrpol\ with a priori bounds and $\mu$ a conformal 
measure for f. Let U be any neighborhood of the
critical point. Then the set

\[ \{ x \in {\bf C} : f^{n} (x) \in U^{c} , \: for \: 
all \: positive \: n \} \]
has zero $\mu$-measure.
\end{lemma}

\proof{
Let us denote the set in the statement of this Lemma by $A$. We have
$A= J(f) \setminus \bigcup_{k} f^{-k}(U)$. So $A$ is a nowhere dense
forward invariant set. Notice that $A \bigcap
\overline {\cal O}$ is empty (because of the definition of $A$ and
because $\overline {\cal O}$ is minimal if $f$ is infinitely many
times renormalizable with a priori bounds). In view of the Lebesgue
density Theorem (see Theorem 2.9.11 in \cite{f}), the set of density
points of $A$ has full measure inside $A$. Here by density points we
mean $x \in A$ such that $lim_{r
\rightarrow 0} {\rm dens}(A|B(x,r))=1$, where $B(x,r)$ is the
Euclidean ball with center at $x$ and radius $r$. Suppose that
$\mu(A)$ is positive. Then we conclude that there exists a density
point $x$ in $A$. There also exists $y$ inside $A$ and a sequence of
natural numbers $k_{j}
\rightarrow \infty$ such that $f^{k_{j}}(x) \rightarrow y$. We can 
pull back a ball of definite size centered in $y$ along $x, f(x), ...,
f^{k_{j}}(x)$ (to be more precise, the size of this ball is ${\rm
dist}(y, \overline {\cal O})$). That implies that we can fix a
positive number $\eta$ and pull back the ball $B(f^{k_{j}}(x),
\eta)$, along $x, f(x), ..., f^{k_{j}}(x)$. Since $A$
is nowhere dense and $\mu$ is positive on non-empty open subsets of
the Julia set, for large $k_{j}$ we have:

\[ \mu (B(f^{k_{j}}(x), \frac{\eta}{2}) \setminus A) \geq  \mu(B(y, 
\frac{\eta}{4}) \setminus A) > 0. \]

As a consequence of Koebe's Theorem, the definition of conformal measure
and the invariance of $A$ we have:

\[ K^{-1} \mu ( B(x, K^{-1} \frac{\eta}{2} |Df^{k_{j}}(x)|^{-1}) \setminus A) 
\leq |Df^{k_{j}}(x)|^{- \delta} \mu ( B(f^{k_{j}}(x),\frac{\eta}{2}) 
\setminus A  ) \leq \] 
\[\leq  K \mu ( B(x, K \frac{\eta}{2} |Df^{k_{j}}(x)|^{-1} ) 
\setminus A). \]

Let us denote $r=K \frac{\eta}{2} |Df^{k_{j}}(x)|^{-1}$. From the
above and from the definition of conformal measure we get:

\[\frac{\mu ( B(x, r ) \setminus A)}{\mu (B(x, r ))} \geq
\frac{|Df^{k_{j}}(x)|^{-\delta} K^{-1} \mu ( B(f^{k_{j}}(x),
\frac{\eta}{2}) 
\setminus A  )}{\mu( B(x, r ))} \geq \]

\[ \geq \frac{K^{-1}|Df^{k_{j}}(x)|^{\delta}}{\mu (f^{k_{j}}(B(x, r )))}
|Df^{k_{j}}(x)|^
{-\delta} K^{-1} \mu ( B(f^{k_{j}}(x),\frac{\eta}{2}) \setminus A )
\geq \]

\[ \geq \frac{ K^{-2} \mu ( B(f^{k_{j}}(x),\frac{\eta}{2}) \setminus A  )}
{\mu (f^{k_{j}}(B(x, r )))}  \geq
\frac{K^{-2} \mu ( B(f^{k_{j}}(x),\frac{ \eta}{2}) \setminus A  )}
{1} \geq\]   

\[\geq K^{-2} \mu(B(y, \frac{\eta}{4}) \setminus A) \geq c > 0 . \]

As $\lim_{k_{j} \rightarrow \infty} |Df^{k_{j}}(x)| = \infty$ 
(because of bounded distortion and lack of normality inside $J(f)$) we get:
\[{\rm limsup}_{r \rightarrow 0} \frac{\mu(B(x,r) \setminus A)}
{\mu(B(x,r))} >0, \]
which contradicts the choice of $x$ as a density point of $A$.
}

\begin{lemma}
\label{zeromeasinfrenormdeglemma}
Let $f(z)=z^{l} +c$, with $l$ even and $c$ real be an 
infinitely renormalizable polynomial and $\mu$ a conformal 
measure for f. Let U be any neighborhood of the
critical point. Then the set

\[ \{ x \in {\bf C} : f^{n} (x) \in U^{c} , \: for \: 
all \: positive \: n \} \]
has zero $\mu$-measure.
\end{lemma}

\proof{
The proof of this Lemma is identical to the proof of the previous
Lemma. The essential information we used in the proof of
Lemma~\ref{zeromeasinfrenormlemma} was the complex bounds. The complex
bounds in the present case is guaranteed by Theorem~\ref{complexboundsthm}. 
}

Note that in Lemma~\ref{zeromeasinfrenormlemma} and
Lemma~\ref{zeromeasinfrenormdeglemma} we used the fact that $f$ restricted to
$\overline {\cal O}$ is minimal which is not necessarily true for polynomials
which are at most finitely many times renormalizable. On the other hand, in
Lemma~\ref{zeromeasyoclemma} we used the fact that we have a partition for the
entire Julia set by puzzle pieces whose pre-images shrink to points. We do not
have that for the polynomials in Lemma~\ref{zeromeasinfrenormlemma} and
Lemma~\ref{zeromeasinfrenormdeglemma}.

From the previous Lemmas we conclude that the set 
\[ {\cal W} = \{ z \in
J(f) : 0 \in w(z) \} \] has full measure, i. e., $\mu ({\cal W} ) =
1$.  Here $w(z)$ denotes the $w$-limit set of $z$. 

Remember that in the case of a finitely many times renormalizable
polynomial (with only repelling periodic points), the first index of the
principal nest is always $0$: $\yocdepzn$. That is because we do not
have renormalization levels.

Let $X \subset {\cal W}$ be any measurable set.  If $f$ is a finitely
many times renormalizable polynomial (real, if the degree is greater
than 2) with only repelling periodic points we can create a cover of
$X$ by puzzle pieces as follows: fix $\yocdepzn$. For every $x \in X$
there exists a first time $m$ such that $f^{m}(x) \in
\yocdepzn$. So we can pull $\yocdepzn$ along the orbit of $x$ back to a puzzle
piece containing $x$. Changing $x \in X$ we will obtain the desired cover. Let
us call this cover $O_{n}$. We can make a similar construction for any
$\seclimb$ polynomial with a priori bounds using the sets $\yocdepmtm$
constructed in Subsection~\ref{lyupol} (we will just consider the
renormalization levels where we have the a priori bounds). For any real
unimodal infinitely renormalizable polynomials we can also repeat the same
construction using the sets $\sulpiecejmu$ introduced on
Theorem~\ref{complexboundsthm}.  We have the following properties:

\begin{romanEnum}

 \item  $O_{n}$ is an open cover;
 \item  $O_{n} \subset O_{n-1}$;
 \item  $\bigcap O_{n} = X$;
 \item $\mu (O_{n}) \rightarrow \mu (X)$ as $n \rightarrow \infty$.

\end{romanEnum}

The first and the second properties are trivial. The third one is a
consequence of Theorem~\ref{yocthm} and Theorem~\ref{lsthm}, if $f$ is
a finitely many times renormalizable polynomial (real or degree two)
with only repelling periodic points. The same fact follows for
$\seclimb$ polynomials with a priori bounds, if we use
Lemma~\ref{unbranchedlemma}. The last one follows by regularity of the
measure $\mu$.

To simplify the notation, elements of $O_{n}$ will be denoted by the
letter $U$ (indexed in some convenient fashion).

\begin{lemma}
\label{existdensptlemma2}
For all i, there exists $U^{i}$ in $O_{i}$ such that \mbox{${\rm
dens}(X | U^{i}) \rightarrow 1$}, as $i \rightarrow \infty$.
\end{lemma}

\proof{
Similar to Lemma~\ref{existdensptlemma}.
}

We will now finish this section with a Lemma for polynomials with
parabolic periodic points:

\begin{lemma}
\label{zeromeasparab}
Let $f(z)=z^{l}+c$, $l$ even and $c$ complex, be a polynomial with a
parabolic periodic point and $\mu$ a conformal measure for f. Let $U$
be any neighborhood of the parabolic periodic point. Then the set

\[ \{ x \in {\bf C} : f^{n} (x) \in U^{c} , \: for \: 
all \: positive \: n \} \]
has zero $\mu$-measure.
\end{lemma}

\proof{We can assume that $f$ has a fixed point $w$. If this is not
the case, we change $f$ by a convenient power $f^{n}$ in order to get
a fixed point.

The orbit of the critical set of $f$ converges to $w$. That is because the
critical set of $f$ (or $f^{n}$) is contained in the union of the
attracting petals of $w$. In that case, if $U$ is any neighborhood of $w$,
there exists an $\varepsilon >0$ such that all the inverse branches of
$f^{m}$, for any natural $m$, is defined in $B(z, \varepsilon)$, for any $z$
outside $U$. If we denote by $A$ the set in the statement of this Lemma, then
the distance from $A$ to the orbit of the critical set is positive.  Now we
can repeat step by step the proof of the Lemma~\ref{zeromeasinfrenormlemma}.
}

\section{Proof of Theorem 1}
\label{prooferg}

\subsection{The non-parabolic cases}

Let $Y \subset {\cal W} = \{ z \in J(f) : 0 \in w(z) \} \subset $
\julf\ be an $f$-invariant set (remember that ${\cal W}$ has full
measure).  Suppose that $\mu (Y) >0$.

If $f$ is a Yoccoz polynomial, by Lemma~\ref{existdensptlemma2}
we can find $U^{n}$ in $O_{n}$ such that ${\rm dens} (Y | U^{n})
\rightarrow 1$. Let $f^{j_(n)} : U^{n} \rightarrow
\yocdepzn$ be an isomorphism (given by the definition of $O_{n}$).
Then by Lemma~\ref{unbranchednonrenlemma} 
and Koebe's Theorem we conclude
that $f^{j_(n)}$ has bounded distortion, i.e.:

\[\frac{1}{K} \leq \frac{|D(f^{j_(n)})(z_{1})|}{|D(f^{j_(n)})(z_{2})|} \leq K\]
for all $z_{1}$ and $z_{2}$ in $U^{n}$, where $K$
depends just on $c(f)$, the constant that appears in the statement of
Lemma~\ref{unbranchednonrenlemma}.

Now let us apply Proposition~\ref{denstransfprop} to the sets $Y^c \bigcap
U^{n}$ and $U^{n}$ with respect to
the map $f^{j_(n)}$. Due to the fact that the set $Y$ is $f$-invariant
and that $f^{j_(n)}(U^{n}) = \yocdepzn$ we get
$\frac{1}{K^ \delta} {\rm dens}(Y^c | U^{n} ) \leq
{\rm dens}(Y^c | \yocdepzn) \leq K^ \delta {\rm dens}(Y^c |
U^{n})$.

We know that ${\rm dens}(Y | U^{n} ) \rightarrow 1$. Passing to the
complement of $Y$ we get ${\rm dens}(Y^c | U^{n} )
\rightarrow 0$. From this and the above inequalities we conclude that
\mbox{${\rm dens}(Y^c | \yocdepzn) \rightarrow 0$}
or ${\rm dens}(Y | \yocdepzn )\rightarrow 1$.

Notice that if $\mu (Y^c) >0$ then we can repeat the argument changing
$Y$ by $Y^c$. Doing this we get ${\rm dens}(Y^c |
\yocdepzn)\rightarrow 1$ and that contradicts the
previous limit because ${\rm dens}(Y|\yocdepzn) + {\rm
dens}(Y^c|\yocdepzn ) =1$.

So we conclude that $\mu (Y^c)$ =0, or equivalently, $\mu (Y) =1$.
This finishes the proof of the Theorem if $f$ is a Yoccoz polynomial.
For any other finitely many times renormalizable real polynomial
with only repelling periodic points, the proof is identical.

If $f$ is a $\seclimb$ polynomial with a priori bounds, the proof of
Theorem~\ref{ergthm} is basically the same. The only difference is
that we use Lemma~\ref{unbranchedlemma} instead of
Lemma~\ref{unbranchednonrenlemma}. If $f(z)=z^{l}+c$ is infinitely many
times renormalizable, $l$ even and $c$ real, then we need
Theorem~\ref{complexboundsthm} and Lemma~\ref{unbranchedinfrenlemma}
to carry out the above argument. Again the proof is the same.

\subsection{The parabolic case}

Now, assume that $f(z)=z^{l}+c$, $l$ even and $c$ complex with a fixed
parabolic point $w$. We know by Lemma~\ref{zeromeasparab} that if $U$
is a neighborhood of $w$, then the set $\cup_{n=0}^{\infty}f^{-n}(U)$
has full measure. Let us assume that $U$ has a small diameter. Then
the set $f^{-1}(U)$ has $l$ connected components: one containing $w$
and the others containing $w_{i}$, the pre-images of $w$, other than
$w$ itself. The connected component of $f^{-1}(U)$ containing $w_{i}$
will be denoted by $U'_{i}$. We will denote the union of all the
$U'_{i}$'s by $U'$. If $U$ is small, then $U \cap J(f)$ is contained
in the union of the repelling petals of the parabolic point $w$. So,
$f^{-1}(U) \subset U
\cup U'$, up to a set of zero measure (remember that $\mu$ is
supported on $J(f)$).

Taking into account all the previous observations, we conclude that up
to a set of zero measure we have:
\[ \cup_{n=0}^{\infty}f^{-n}(U) = \cup_{n=0}^{\infty}f^{-n}(U') \cup U \]

\noindent{\bf $1^{st}$ case}: Let us assume that $\mu(w)=0$. Then by
regularity of $\mu$ and by the last equality we have that the following is
true: for any $\varepsilon >0$, there exists a sufficiently small neighborhood
$U$ of $w$ such that $\mu(\cup_{n=0}^{\infty}f^{-n}(U')) > 1 -\varepsilon$
(just take $U$ such that $\mu (U) < \varepsilon$).

If $U_{2} \subset U_{1}$ then
$\cup_{n=0}^{\infty}f^{-n}(U_{2}')
\subset \cup_{n=0}^{\infty}f^{-n}(U_{1}')$. So 
$\mu(\cup_{n=0}^{\infty}f^{-n}(U_{1}')) \geq
\mu(\cup_{n=0}^{\infty}f^{-n}(U_{2}')) > 1-\varepsilon_{2}$, where
$\varepsilon_{2}$ depends on $U_{2}$. So if $\diam (U_{2})$ is taken
arbitrarily small, then $\varepsilon_{2}$ will be arbitrarily close to
zero.
Then we conclude that $\mu(\cup_{n=0}^{\infty}f^{-n}(U_{1}'))=1$.

As $f^{n}(0) \rightarrow w$, the points $w_{i}$ are at a positive
distance from the critical orbit. This implies that there exists a
positive number $\alpha$ such that all the inverse branches of $f^{n}$
are defined on the ball $B(w_{i},
\alpha)$.

Let us show that $\mu$ is an ergodic measure. Let $Y \subset J(f)$ be
an $f$-invariant set such that $\mu (Y) >0$. In view of the Lebesgue
density Theorem (see Theorem 2.9.11 in \cite{f}) the set of density
points of $Y$ has full measure inside $Y$. 

By the previous paragraphs we conclude that there exist $x \in Y$ and
a sequence $\{ k_{j} \}$ of numbers such that $\lim_{r \rightarrow 0}
{\rm dens} (Y|B(x,r))=1$ and that $\lim_{k_{j} \rightarrow \infty}
f^{k_{j}}(x) = w_{i_{0}}$, for some fixed $i_{0}$.

We will show that $w_{i_{0}}$ is a density point of $Y$. If this is not the
case, then $\lim_{r_{i} \rightarrow 0} {\rm den}s (Y|B(-w, r_{i}))< 1$, for
some sequence of positive numbers $r_{i}$ tending to zero. So there
exists $\eta < \alpha$ such that $\mu (B(w_{i_{0}}, \frac{\eta}{4}) \setminus
Y)>0$. Now we can mimic the proof of
Lemma~\ref{zeromeasinfrenormlemma} to conclude that $x$ is not a
density point of $Y$. Contradiction! So we conclude that $w_{i_{0}}$ is a
density point of $Y$.

As the restriction of $f$ to a small neighborhood of $w_{i_{0}}$ is an
isomorphism onto a neighborhood of $w$ and $Y$ is $f$-invariant, we
conclude that $w$ is a density point of $Y$.

If the measure of the complement of $Y$ is positive, then we can show
that $w$ is a density point of the complement of $Y$. We do that just
changing $Y$ by its complement in the above reasoning. So we must have
$\mu(Y)=1$, and $\mu$ is ergodic.

\noindent{\bf $2^{nd}$ case}: If $\mu (w)>0$, then Lemma 11 and Theorem 13
from \cite{du2} shows that $\mu$ is supported on the grand-orbit of $w$, and
then it is ergodic.

If $f$ has a parabolic periodic point $w$ of period $p$, then $f^{p}$
has a fixed parabolic point. We can prove the Theorem in exactly the
same way we did before.
So we proved Theorem~\ref{ergthm}.

\bibliographystyle{alpha}

\end{document}